\documentclass{amsart}
\usepackage{amsmath}
\usepackage{amscd}
\usepackage{amssymb}

\theoremstyle{theorem}
\newtheorem{theorem}{Theorem}

\theoremstyle{corollary}
\newtheorem*{corollary}{Corollary}
\theoremstyle{definition}
\newtheorem{definition}{Definition}
\newtheorem*{lemma}{Lemma}
\newtheorem*{remark}{Remark}
\newtheorem*{example}{Example}
\newtheorem*{affil}{}
\def\men{\smallsetminus}

\begin{document}

\title{On punctured locally compact spaces}

\author{Giuseppe De Marco}

\maketitle

\begin{abstract}
In a recent paper \cite{T} the fact that  a class of locally compact metric spaces $X$, among which are Euclidean spaces,  are not homemorphic to their punctured version $X\men\{p\}$, was given an interesting new proof which does not use algebraic topology; essential tools of this proof are  a boundedly compact metric structure, 
and path--connectedness near infinity. Here we show that local compactness and ordinary connectedness near infinity suffice; no metrizability is needed, and moreover we can also delete whole compact subsets, not only single points.  Some non--homeomorphism results on many--holed Euclidean balls are also obtained. This note ought to distil the essence of the method developed in \cite{T}.
\end{abstract}
\noindent

\section*{One--point compactification.} This notion is central in the note, so we expand it a bit. The {\em one--point compactification}, or Alexandroff compactification, of a topological space $X$ is the space $\alpha X=X\cup\{\infty_X\}$ ($\infty_X$, the {\em point at infinity}, is simply an object not in $X$) whose topology consists of the open sets of $X$, and all sets of the form $\alpha X\men K$, where $K$ is a closed compact subset of $X$. If $X$ is compact  then $\infty_X$ is an isolated point of $\alpha X$; otherwise $\alpha X$ is a compact space that has $X$ as an open dense subspace. 

The one--point compactification is familiar to analysts and topologists alike: a classical fact is that $\alpha \mathbb R^n$ is homeomorphic to the $n-$sphere $\mathbb S^n=\{x\in\mathbb R^{n+1}:\,|x|=1\}$, via stereographic projection on $\mathbb R^n$ from the north pole of $\mathbb S^n$ (\cite{B}, 4.3). 

Notice that a neighbourhood (nbhd) of the point at infinity is a set of the form $\alpha X\men A$=$(X\men A)\cup\{\infty_X\}$, where $A$ is a relatively compact subset of $X$, i.e. such that $\operatorname{cl}_X(A)$  is compact. Even if the space $X$ is  separated (=Hausdorff, or $T_2$) the space  $\alpha X$, although compact, need not be separated. For instance, if $X$ is an infinite dimensional normed space, then no $x\in X$ has a disjoint nbhd with $\infty_X$: recall that by Riesz's theorem compact subsets in an infinite dimensional normed space have empty interior (\cite{D}, 5.5.9). In fact, the points of a space $X$ which have a nbhd disjoint from a nbhd of $\infty_X$ are exactly those which have a closed compact nbhd in $X$, as is immediate by the definition. Thus to make $\alpha X$ separated  we need $X$  separated (of course!) but also that every point of $X$ has in $X$ a compact nbhd (in separated spaces compact subsets are closed, so we don't need to add ''closed''). This condition is called {\em local compactness}. It can be proved that  in a locally compact separated space the compact nbhds of a point are a base for the nbhd system of the point (\cite{K}, 5.17). This implies that open subspaces of separated locally compact spaces are locally compact. An easy compactness argument shows also that every compact subset of a separated locally compact space has a nbhd base consisting of compact sets (\cite{K}, 5.18).
 
  In all the paper, ''locally compact and separated'' could be replaced by "locally compact and regular'': what is needed is that at every point  there is  a base for the nbhd system consisting of closed compact sets. But the Hausdorff axiom is much more familiar, and spares us the pain of saying ''closed compact''  instead of simply ''compact''.
  
  Hereafter, locally compact spaces are  also assumed separated. 

\section*{Main result.} If $U$ is a nbhd of a point $p$ in a space $X$, then $U\men\{p\}$ is the associated {\em punctured} nbhd (called  {\em deleted} nbhd  in some limit definitions).
\begin{lemma} Let $X$ be a  locally compact non--compact space, let $C$ be a compact subset of $X$, and let $Y=X\men C$. Then there exist in $X$ open disjoint subsets $U_0,\, V_0$, with $C\subseteq U_0$ and $X\men V_0$  compact, such that every punctured nbhd of $\infty_Y$ contains a punctured nbhd of the form
 \[\tag{*}(U\men C)\cup V\quad\text{with}\quad U\subseteq U_0,\, V\subseteq V_0,\]
 where $U$ is a nbhd of $C$ in $X$ and $V$ a punctured nbhd of $\infty_X$. And every set of this form is a punctured nbhd of $\infty_Y$.\end{lemma}
 Before the proof: the lemma says that removing a compact piece $C$ from $X$, the new point at infinity $\infty_Y$ will have punctured nbhds consisting of two parts, one around the old point at infinity $\infty_X$, and another around the excised part $C$.
 \begin{proof} Since $Y$ is open in $X$, $Y$ is locally compact. Let $W_0$ be a compact  nbhd  of $C$ in $X$, and let $B_0$ be the boundary of $W_0$ in $X$. Then $B_0$, a closed subset of the compact set $W_0$, is compact, and $X\men B_0$ is the disjoint union of the open sets $U_0$ and  $V_0$, where $U_0=\operatorname{int}_XW_0$, and $V_0=Y\men W_0$; notice that $C\subseteq U_0$. Then $P_0=(U_0\men C)\cup V_0=Y\men B_0$ is a punctured nbhd of $\infty_Y$. If $P$ is another punctured nbhd of $\infty_Y$, $P$ contains a punctured open nbhd of the form $Q=Y\men K$, where $K$ is a compact subset of $Y$ containing $B_0$, so that $Q\subseteq (U_0\men C)\cup V_0$; then $U=(Q\cup C)\cap U_0$ and $V=Q\cap V_0$ are as required. The last observation is easy: it is not restrictive to assume that $U$ and $V$ are open, and 
 \[Y\men ((U\men C)\cup V)=(Y\cup C)\men((U\men C)\cup V\cup C)=X\men(U\cup V)\]
 is compact.
 \end{proof} 
 \begin{remark} If $C$ is not open, then $U\men C$ in the lemma above cannot be empty. For  $U\men C=\emptyset$ means that $U=C$, i.e. $C$ is a nbhd of itself, i.e. open. And no $V$ can be empty, since $X$ is non--compact. Thus no punctured nbhd of the form (*) can be connected.\end{remark}
 \begin{definition} We say that a topological space $X$ is {\em connected near infinity} if for every relatively compact subset $A$ of $X$ there is a relatively compact subset $B$ of $X$ containing $A$ such that $X\men B$ is connected.\end{definition}
 This condition is vacuously satisfied for $X$ compact (the empty space is connected!). In the non--compact case, it means that  every punctured nbhd of $\infty_X$ contains a connected one. Since homeomorphisms of  spaces trivially  extend to homeomorphisms of their one point compactifications, by mapping  points at infinity into each other, a space homeomorphic to a space connected near infinity is also connected near infinity. Typical spaces connected near infinity are normed linear spaces of real dimension $>1$; as remarked in the previous section, only the finite dimensional ones have a separated one--point compactification. A locally compact non--compact non--metrizable space connected and connected  near infinity is the long half--line, $\omega_1\times[0,1[$, lexicographically ordered, with the order topology (see \cite{K}, 5K; $\omega_1$ is the first uncountable ordinal).
 
 \begin{theorem} Let $X$ be a  locally compact non--compact space, and let $C\subseteq X$ be a compact non--open subset of $X$. Then $Y=X\men C$ is not connected near infinity.\end{theorem}
  \begin{proof}  Apply the lemma above, keeping account of the Remark.\end{proof}
    \begin{corollary} Let $X$ be a locally compact space connected near infinity, and let $C$ be compact non--open in $X$. Then $X\men C$ is not homemorphic to $X$.  In particular, if $p\in X$ is non--isolated in $X$, then $X\men\{p\}$ is not homemorphic to $X$.\end{corollary}
 \begin{proof} Since $C$ is not open in $X$, $X\men C$ is not compact, so that it is not homeomorphic to $X$ if $X$ is compact. If $X$ is non--compact,  Theorem 1 says that $X\men C$ is not connected near infinity, hence it is not homeomorphic to $X$. For the last assertion, simply recall that  $p$   isolated in $X$ means that $\{p\}$ is open in $X$. \end{proof}
  Hence  if $n\ge 2$, and $K$ is a non--empty compact subset of $\mathbb R^n$, then $\mathbb R^n\men K$ is not homeomorphic to $\mathbb R^n$.
 
 \begin{example} If $p$  is isolated in $X$, then easy examples show that $X$ may be connected near infinity, and $X\men\{p\}$  homeomorphic to $X$. Take e.g. the subspace of the reals  $X=\{x_n:n\in\mathbb{N}\}\cup[1,\infty[$ with  $x_n=n/(n+1)$; then $x_0=0$ is isolated, and the map $\phi: X\men\{0\}\to X$ given by $\phi(x_n)=x_{n-1}$, $\phi(x)=x$ for $x\in[1,\infty[$ is plainly a homeomorphism. But there is no homeomorphism between  $X$  and $X\men\{1\}$; this follows from the above corollary. There is also  this alternative proof: a homeomorphism  maps isolated points into isolated points,  so that  a homeomorphism $\psi:X\to X \men\{1\}$  induces a permutation of $\{x_n:n\in\mathbb{N}\}$; then $\psi(1)=\lim_n\psi(x_n)=\lim_mx_m=1$ (a rearrangement of a converging sequence converges to the same limit)  contradicting $1\notin X\men\{1\}$.
  \end{example}
  \section*{Topological sums.}In the last two sections  {\em topological sums} appear. A space $X$ admits a splitting into a {\em topological sum} of a family $(X_\lambda)_{\lambda \in\Lambda}$ of subspaces  if $X=\biguplus_{\lambda \in\Lambda}X_\lambda$, disjoint union, where all $X_\lambda$ are clopen (=open and closed) subspaces of $X$ (see \cite{B}, 8.11); call the splitting non--trivial if all spaces $X_\lambda$ are non--empty, and $\Lambda$ is not a singleton. By definition, a space is connected if it has no non--trivial splitting.  Observe now that if $K$ is a compact subset of $X$, then there is a finite subset $F$ of $\Lambda$ such that $K\subseteq\bigcup_{\lambda\in F}X_\lambda$, and $K_\lambda=K\cap X_\lambda$ is compact for every $\lambda\in\Lambda$. Then every punctured nbhd of $\infty_X$ contains a set of the form 
\[X\men K=\left(\biguplus_{\lambda\in F}U_\lambda\right)\uplus\left(\biguplus_{\lambda\in\Lambda\men F}X_\lambda\right),\]
where $U_\lambda=X_\lambda\men K$ is a punctured nbhd of $\infty_{X_\lambda}$, for every $\lambda\in F$.
Hence  $X$ is not connected near infinity if $\Lambda$ is infinite, or if more than one of the $X_\lambda$'s is non--compact.

A connected non--empty subset of $X$ is contained in exactly one $X_\lambda$.

 \section*{A generalization.} Take the open unit disk in the plane, and consider the spaces obtained by removing from it $p$ points, respectively $q$ points, where $p,q\in\mathbb N$. Are these spaces homeomorphic if $p\ne q$? Algebraic topology says no, because their fundamental groups are non--isomorphic. We see here that the previous ideas can be adapted to give the same answer; and our method is (in this case) even more sensitive, because it can distinguish between the open holed disk and the closed holed disk, unlike the fundamental group!

  \begin{definition} A topological space $X$ is said to be $m-$disconnected near infinity if every punctured nbhd of $\infty_X$ contains a punctured nbhd which has exactly $m+1$ connected components, and $m$ is the smallest finite cardinal for which this happens.\end{definition}
  Thus $0-$disconnected near infinity means connected near infinity, the real line is $1-$disconnected near infinity, the union of coordinate axes in $\mathbb R^2$ is $3-$disconnected near infinity, etc. (and we might say that compact spaces are $-1$ disconnected near infinity!). It is clear that  this property is a topological invariant, i.e. it is shared by homeomorphic spaces.
  
  Using the  Lemma above it is easy to prove the following:
   \begin{theorem} Let $X$ be a  locally compact  space, and let $C\subseteq X$ be a compact non--open subset of $X$. Assume that  $X$ is $m$--disconnected near infinity, for some natural number $m$. Then $Y=X\men C$ is not $m$--disconnected near infinity, and therefore not homeomorphic to $X$.\end{theorem}
   Of course in general $Y$ will not be $(m+1)-$disconnected  near infinity: this depends on how the ''hole'' $C$ is placed in $X$.

  If a locally compact space $X$ has an open relatively compact subset $A$ such that   $X\men A$ is the topological sum of $m$ non--compact spaces $X_1,\dots,X_m$, all connected near infinity, then $X$ is  $(m-1)-$disconnected near infinity; in fact it is easy, arguing as in the proof of lemma 1, to see that for every punctured nbhd $U$ of $\infty_X$ contained in $X\men A$ we have $U=\biguplus_{j=1}^mU_j$, where each $U_j$ is a punctured nbhd of $\infty_{X_j}$. This could be used to prove the following theorem, but the simple direct argument allowed by the metric structure is perhaps more convincing.  
  
  \begin{theorem} Let $D$ be   the open unit ball in $\mathbb R^n$, $n\ge2$, and consider the space $X(p)=D\men\{a_1,\dots, a_p\}$, where $p$ is a finite cardinal. Then $X(p)$ is $p$--disconnected near infinity; and $Y(p)=\bar D\men\{a_1,\dots, a_p\}$ is $(p-1)$--disconnected near infinity.\end{theorem}
  \begin{proof} Denote by $\rho$ the distance function from  the complement of $X(p)$, that is $\rho(x)=\inf\{|x-y|:\,y\in\mathbb R^n\men X(p)\}$. If $K$ is a compact subset of $X(p)$, $\rho$ has a strictly positive minimum $\delta$ on $K$, so every nbhd of $\infty_{X(p)}$ contains a set  such that $X_\delta=\{x\in X(p):\,\rho(x)<\delta\}(\subseteq X(p)\men K)$. If $\delta>0$ is small enough, to be precise smaller than $\min\{(1-|a_j|)/2, |a_j-a_k|/2,\,j\ne k, j,\,k= 1,\dots,p\}$ then $X_\delta$ is the union of $p+1$ pairwise disjoint open sets:
  \[ X_\delta=
  \bigcup_{j=1}^p(B_\delta(a_j)\men\{a_j\})\cup\{x\in\mathbb R^n:\,1-\delta<|x|<1\}.\]
  Since $n\ge2$ these sets are all connected, so that they are the connected components of $X_\delta$. Now a punctured nbhd $V$ of $\infty_{X(p)}$ contained in $X_\delta$ contains $X_\sigma$ for some $\sigma<\delta$, hence it intersects every component of $X_\delta$, so that it cannot have less than $p+1$ connected components. The proof for $Y(p)$ is analogous; punctured nbhds of $\infty_{Y(p)}$ lack the ''spherical shell'' at the boundary of $D$, so  they have $p$ components.  
  \end{proof}

  
  \section*{Comments on paper \cite{T}.} The idea for the present note came from that paper.  The main theorem there is stated (verbatim) as follows:
  
  {\sl Let $(X,\,d)$ be a complete, locally compact metric space and assume that $X\men B_R(q)$ is path connected, for some $q\in X$ and every $R>0$ large enough. Then there is no homeomorphism between $X$ and its punctured version $X\men\{p\}$.}

  The proof given in \cite{T} assumes that every closed bounded subset of $X$ is compact. The very first argument in the proof in fact is: if $y_j$ is a sequence in $X$ such that $d(p,\,y_j)$ remains bounded, then $y_j$ has a converging subsequence. Metric spaces with this property are called by some {\em boundedly compact}. The only metric requirement made above on the topologizing metric  is completeness, and even coupled with local compactness   it cannot imply bounded compactness. For we know that if $d$ is a metric on $X$, then also $d_1(x,\,y)=\min\{d(x,\,y),1\}$ is a metric on $X$ with the same topology  and the same Cauchy sequences, but clearly $X$, although now of $d_1-$diameter $1$ or less, has not become compact! Boundedly compact spaces are always complete  and  locally compact;  what is really proved in \cite{T} is (correcting also the minor oversight on isolated points):
  
  {\sl Let $(X,\,d)$ be a boundedly  compact metric space and assume that $X\men B_R(q)$ is path connected, for some $q\in X$ and every $R>0$ large enough. Then there is no homeomorphism between $X$ and its punctured version $X\men\{p\}$, if $p$ is non--isolated in $X$.}

Of course the real question is: does a complete  locally compact  metric space $X$   admit a topologically equivalent boundedly compact metric? The answer is negative: a boundedly compact metric space is $\sigma$--compact, i.e. a countable union of compact subspaces, since  $X=\bigcup_{n=1}^\infty\operatorname{cl}_X(B_n(q))$, with $q$ any point of $X$; then it is a Lindel\"of space, and hence it has a countable base for the topology; see \cite{K} for all these notions and implications. Since all subspaces of $X$ have then a countable base, no space that contains a discrete uncountable subspace can admit a boundedly compact metric. In particular,  uncountable discrete spaces, although locally compact, and complete under the usual $\{0,\,1\}$--valued discrete metric, do not admit such a metric. It is interesting to note that a locally compact space has a metrizable one-point compactification if and only if its topology is induced by a boundedly compact metric.  This fact is certainly well--known, as are all the preceding ones,  but I cannot pinpoint a satisfactory reference, so here is a proof. If $X\cup\{\infty_X\}$ is topologized by a metric $\rho$, it is easy to see that the metric $d:X\times X\to [0,\infty[$ given by:
\[d(x,\,y)=\rho(x,\,y)+\left|\dfrac1{\rho(x,\infty_X)}-\dfrac1{\rho(y,\infty_X)}\right|\]
is boundedly compact and topologizes $X$. And if $X$ is locally compact and second countable, then so  trivially is $\alpha X$ (\cite{B}, 12.12)  which is then metrizable,  by Urysohn metrization theorem.

But this is not yet the whole story: a paracompact  locally compact space $X$ can be split into a topological sum of (non--empty) $\sigma$--compact subspaces, $X=\biguplus_{\lambda \in\Lambda}X_\lambda$, with each $X_\lambda$ $\sigma$--compact (\cite{Bo}, Th\'eoreme 5). Since metrizable spaces are paracompact, if $X$ is metrizable all the $X_\lambda$ admit a boundedly compact metric; and any connected subset of $X$ will be contained in a single $X_\lambda$. If $X$ is connected, the splitting is trivial: connected locally compact metrizable spaces are $\sigma$--compact, and admit a boundedly compact metric. And as observed previously,  if there is more than one non--compact $X_\lambda$, or if $\Lambda$ is infinite, the space will not be connected near infinity. 

It was the desire to understand if boundedly compact metrizability was essential for the validity of the result that spawned the present paper.

\begin{affil}{\small
Dipartimento di Matematica ''Tullio Levi--Civita'', Universit\`a di Padova, Via Trieste 63, 35121 Padova, Italy\\
gdemarco@math.unipd.it}
\end{affil}

MSC: 54-02

\end{document}